\documentclass
{article}

\usepackage{textcomp}
\usepackage{graphicx}

\usepackage{latexsym}
\usepackage{amsfonts}
\usepackage[psamsfonts]{amssymb}
\usepackage{enumerate}

\RequirePackage[OT1]{fontenc}
\RequirePackage[amsthm,amsmath]{imsart}

\theoremstyle{plain} 
\newtheorem{theorem}{Theorem}
\newtheorem{corollary}{Corollary}
\newtheorem{lemma}{Lemma}
{Proposition}

\theoremstyle{definition} 

\theoremstyle{definition} 

\theoremstyle{remark} 

\theoremstyle{remark} 
{Remark}
\newtheorem*{remark*}{Remark}

\renewcommand{\u}{\nearrow}
\renewcommand{\d}{\searrow}

\newcommand{\lp}{\left(}
\newcommand{\rp}{\right)}

\newcommand{\p}[2]{\frac{\partial #1}{\partial #2}}

\renewcommand{\le}{\leqslant}
\renewcommand{\ge}{\geqslant}

\newcommand{\la}{\lambda}
\newcommand{\de}{\delta}

\newcommand{\vpi}{\varphi}
\newcommand{\ps}{\overline{\Phi}}

\renewcommand{\P}{\operatorname{\mathsf{P}}} 

\newcommand{\E}{\operatorname{\mathsf{E}}}

\newcommand{\LLe}{\operatorname{\mathrm{LL\!e}}}
\newcommand{\LG}{\operatorname{\mathrm{LG}}}
\newcommand{\GL}{\operatorname{\mathrm{GL}}}
\newcommand{\GG}[1]{\operatorname{\mathrm{GG}_{#1}}}
\newcommand{\GE}{\operatorname{\mathrm{GE}}}
\newcommand{\ELe}{\operatorname{\mathrm{EL\!e}}}
\newcommand{\EG}[1]{\operatorname{\mathrm{EG}_{#1}}}
\newcommand{\A}[1]{\operatorname{\mathrm{A}_{#1}}}
\newcommand{\X}[1]{\operatorname{\mathrm{X}_{#1}}}

\newcommand{\R}{\mathbb{R}}

\newcommand{\vp}{\varepsilon}

\begin{document}


\begin{frontmatter}

\title{Toward the best constant factor for the Rademacher-Gaussian tail comparison}
\runtitle{Rademacher-Gaussian tail comparison}

\begin{aug}
\author{\fnms{Iosif} \snm{Pinelis}\ead[label=e1]{ipinelis@math.mtu.edu}}
\runauthor{Iosif Pinelis}


\address{Department of Mathematical Sciences\\
Michigan Technological University\\
Houghton, Michigan 49931, USA\\
E-mail: \printead[ipinelis@mtu.edu]{e1}}
\end{aug}






\begin{abstract}
Let $S_n:=a_1\vp_1+\cdots+a_n\vp_n$, where 
$\vp_1,\dots,\vp_n$ are independent Rademacher random variables and 
$a_1,\dots,a_n$ are any real numbers such that
$a_1^2+\dots+a_n^2=1$.
Let $Z\sim N(0,1)$.
It is proved that the best constant factor $c$ in 
inequality
	$\P(S_n\ge x)\le c\P(Z\ge x)$ $\forall x\in\R$ 
	is between two explicitly defined absolute constants $c_1$ and $c_2$ such that $c_1<c_2\approx1.01\,c_1$.
\end{abstract}







\end{frontmatter}


\section{Introduction and summary}\label{intro}

Let $\vp_1,\dots,\vp_n$ be independent Rademacher random variables (r.v.'s), so that $\P(\vp_i=1)=\P(\vp_i=-1)=\frac12$ for all $i$. Let 
$$S_n:=a_1\vp_1+\cdots+a_n\vp_n,$$
where
$a_1,\dots,a_n$ are any real numbers such that
\begin{equation*}
	a_1^2+\dots+a_n^2=1.
\end{equation*}

The best upper exponential bound on the tail probability $\P(Z\ge x)$ for a standard normal random variable $Z$ and a nonnegative number $x$ is $\inf_{t\ge0}e^{-tx}\E e^{tZ}=e^{-x^2/2}$. Thus, a factor of the order of magnitude of $\frac1x$ is ``missing" in this bound, compared with the asymptotics $\P(Z\ge x)\sim\frac1x\,\vpi(x)$ as $x\to\infty$, where $\vpi(x):=e^{-x^2/2}/\sqrt{2\pi}$ is the density function of $Z$. Now it should be clear that any exponential upper bound on the tail probabilities for sums of independent random variables must be missing the $\frac1x$ factor. 

Eaton~\cite{eaton} obtained an upper bound on $\P(S_n\ge x)$, which is asymptotic to $c_3\P(Z\ge x)$ as $x\to\infty$, where 
$$c_3:=\frac{2e^3}9\approx4.46,$$
and he conjectured that $\P(S_n\ge x)\le c_3\frac1x\,\vpi(x)$ for $x>\sqrt2$. 
The stronger form of this conjecture, 
\begin{equation}\label{eq:pin94}
	\P(S_n\ge x)\le c\P(Z\ge x)
\end{equation}
for all $x\in\R$ with $c=c_3$
was proved by Pinelis~\cite{pin94}, along with a multidimensional extension. 

Edelman \cite{edel} proposed inequality
$\P(S_n\ge x) \le \P\lp Z\ge x-1.5/x\rp$ for all $x>0$, but his proof appears to have a gap. A more precise upper bound, with $\ln c_3=1.495\dots$ in place of $1.5$, was recently shown \cite{pin-ed} to be a rather easy corollary of the mentioned result of \cite{pin94}.
Various generalizations and improvements of inequality \eqref{eq:pin94} as well as related results were given by Pinelis \cite{pin98,pin99,pin-eaton,binom,normal,asymm,pin-ed}
and Bentkus \cite{bent-liet02,bent-jtp,bent-ap}.

Bobkov, G\"{o}tze and Houdr\'{e} (BGH) \cite{bgh} gave a simple proof of \eqref{eq:pin94} with a constant factor $c\approx12.01$. Their method was based on the Chapman-Kolmogorov identity for the Markov chain $(S_n)$. Such an identity was used, e.g., in \cite{maximal} to disprove a conjecture by Graversen and Pe\v skir \cite{g-peskir} on $\max_{k\le n}|S_k|$. 

In this paper, we shall show that a modification of the BGH method can be used to obtain inequality \eqref{eq:pin94} with a constant factor $c\approx1.01\,c_*$, where $c_*$ is the best (that is, the smallest) possible constant factor $c$ in \eqref{eq:pin94}.


Let 
$\ps$ and $r$ denote the 
tail function of 
$Z$ and the inverse Mills ratio: 
$$
\text{
$\ps(x):=\P(Z\ge x)=\int_x^\infty\vpi(u)\,du$
\quad and\quad
$r:=\frac\vpi{\rule{3pt}{0pt} \rule{0pt}{10pt} \ps\ }$.
}
$$

\begin{theorem}[Main] \label{th:main}
For the least possible absolute constant factor $c_*$ in inequality \eqref{eq:pin94} one has
\begin{gather*}
c_*\in[c_1,c_2]\approx[3.18,3.22],\quad\text{where}
\\	
c_1:=\frac1{4\ps(\sqrt2)}\quad\text{and}\quad
c_2:=c_1\cdot\Big(1+\tfrac1{250}\big(1+r(\sqrt3\,)\big)\Big)
	\approx c_1\cdot1.01. 
\end{gather*}
\end{theorem}

Here we shall present just one application of Theorem~\ref{th:main}, to self-normalized sums
\begin{equation*}
V_n:=\frac{X_1+\dots+X_n}{\sqrt{X_1^2+\dots+X_n^2}},
\end{equation*}
where, following 
Efron \cite{efron}, we assume that the $X_i$'s satisfy the so-called orthant symmetry condition:
the joint distribution of $\de_1X_1,\dots,\de_n X_n$ is the same for any choice of signs $\de_1,\dots,\de_n\in\{1,-1\}$, so that, in particular, each $X_i$ is symmetrically distributed. It suffices that the $X_i$'s be independent and symmetrically (but not necessarily identically) distributed. 
In particular, $V_n=S_n$ if $X_i=a_i\vp_i$ $\forall i$.
It was noted by Efron that (i) Student's statistic $T_n$ is a monotonic function of the so-called self-normalized sum: $T_n=\sqrt{\frac{n-1}n}\,V_n/\sqrt{1-V_n^2/n}$ and (ii)
the orthant symmetry implies in general that the distribution of $V_n$ is a mixture of the distributions of normalized Rademacher sums $S_n$. Thus, one obtains

\begin{corollary}\label{cor:}
Theorem~\ref{th:main} holds with $V_n$ in place of $S_n$.
\end{corollary}

\section{Proof of Theorem~\ref{th:main}}
\label{proof}

Theorem~\ref{th:main} follows immediately from
Lemma~\ref{lem:c1}, Theorem~\ref{th:refined}, and Lemma~\ref{lem:h1<h}, stated in Subsection~\ref{lemmas, proof of th} below. In particular, Lemma~\ref{lem:h1<h} implies that the upper bound $h_1(x)$ on $\P(S_n\ge x)$ provided by Theorem~\ref{th:refined} 
is somewhat better than the upper bound $c_2\P(Z\ge x)$, implied by Theorem~\ref{th:main}. 

\begin{figure}[htb]
  \begin{center}
{\includegraphics{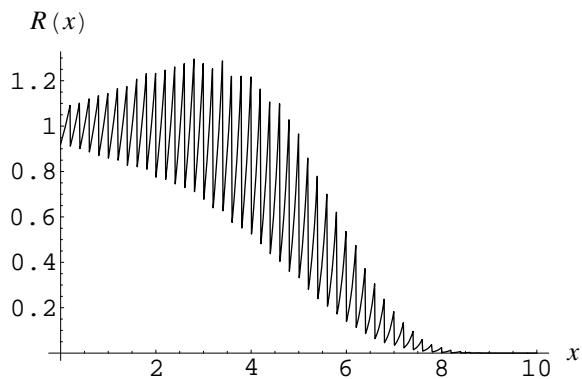} }
  \end{center}
  \caption{
  Ratio of the Rademacher and Gaussian tails for $n=100$ and $a_1=\dots=a_{100}=\frac1{10}$.
}
\label{fig:R(x)}
\end{figure}

While $S_n$ represents a simplest case of the sum of independent non-identically distributed r.v.'s, it is still very difficult to control in a precise manner. Figure~\ref{fig:R(x)} shows the graph of the ratio $R(x):=\P(S_n\ge x)/\P(Z\ge x)$ for $n=100$ and $a_1=\dots=a_n$. 

One can see that even for such a fairly large value of $n$ and equal coefficients $a_1,\dots,a_n$, ratio $R(x)$ oscillates rather wildly. In view of this, the existence of a high-precision inductive argument in the general setting with possibly unequal $a_i$'s may seem very unlikely.
However, such an argument will be presented in this paper.
A key idea in the proof of Theorem~\ref{th:main} is the construction of the upper bound $h_1$ and, in particular, the function $g$ in \eqref{eq:g} and \eqref{eq:h1}, which allows an inductive argument to prove Theorem~\ref{th:refined}. 

The proof of Theorem~\ref{th:refined} is based on a number of lemmas. 
The proofs of all lemmas are deferred to Subsection~\ref{proofs of lemmas}.

\subsection{Statements of lemmas and the proof of Theorem~\ref{th:refined}}
\label{lemmas, proof of th} 

\begin{lemma}\label{lem:c1}
One has $c_*\ge c_1$.
\end{lemma}

For $a\in[0,1)$ and $x\in\R$, introduce
\begin{align}
g(x)&:=c_1\cdot\Big(1+\tfrac1{250}\big(1+r(x)\big)\Big)\ps(x)
=\tfrac{c_1}{250}\cdot\big(251\,\ps(x)+\vpi(x)\big);
\label{eq:g}\\
h(x)&:=c_1\cdot\Big(1+\tfrac1{250}\big(1+r(\sqrt3\,)\big)\Big)\ps(x)=c_2\cdot\ps(x);\notag\\
h_1(x)&:=
	\begin{cases}
	1	&	\text{if }x\le0,\\
	\frac12	&	\text{if }0<x\le1,\\
	\frac1{2x^2}	&	\text{if }1\le x<\sqrt2,\\
	g(x)	&	\text{if }\sqrt2\le x\le\sqrt3,\\
	h(x)	&	\text{if }x\ge\sqrt3;
	\end{cases} 
	\label{eq:h1}\\
K(a,x)&:=h_1(u)+h_1(v)-2h_1(x),\quad\text{where}\notag
\\
u&:=u(a,x):=\frac{x-a}{\sqrt{1-a^2}}\quad\text{and}\label{eq:u}\\
v&:=v(a,x):=\frac{x+a}{\sqrt{1-a^2}}.	\label{eq:v}
\end{align}

\begin{theorem}[Refined]\label{th:refined}
One has 
\begin{equation}\label{eq:le h_1(x)}
\P(S_n\ge x)\le h_1(x)	
\end{equation}
for all $x\in\R$.
\end{theorem}

\begin{lemma}\label{lem:g?h}
One has $g\le h$ on $(-\infty,\sqrt3\,]$ and $g\ge h$ on $[\sqrt3,\infty)$.
\end{lemma}

\begin{lemma}\label{lem:h1<h}
One has $h_1\le h$ on $\R$.
\end{lemma}

\begin{lemma}\label{lem:sqrt3>x>sqrt2}
One has $K(a,x)\le0$ for all $(a,x)\in[0,1)\times[\sqrt2,\sqrt3\,]$.
\end{lemma}

\begin{lemma}\label{lem:x>sqrt3}
One has $K(a,x)\le0$ for all $(a,x)\in[0,1)\times[\sqrt3,\infty)$.
\end{lemma}

Now we can present

\begin{proof}[Proof of Theorem~\ref{th:refined}]
Theorem~\ref{th:refined} will be proved by induction in $n$. 
It is obvious for $n=1$. 
Let now $n\in\{2,3,\dots\}$ and assume that Theorem~\ref{th:refined} holds with $n-1$ in place of $n$.

Note that for $x\le0$ inequality \eqref{eq:le h_1(x)} is trivial.
For $x\in(0,\sqrt2)$, it follows by the symmetry of $S_n$ and Chebyshev's inequality. 
Therefore, assume without loss of generality that 
$x\ge\sqrt2$
and $0\le a_n<1$.
By the Chapman-Kolmogorov identity and induction, 
\begin{align*}
	\P(S_n\ge x)
	&=\tfrac12\,\P(S_{n-1}\ge x-a_n)+\tfrac12\,\P(S_{n-1}\ge x+a_n) \\	&\le\tfrac12\,h_1\big(u(a_n,x)\big)+\tfrac12\,h_1\big(v(a_n,x)\big) \\
	&= h_1(x)+\tfrac12\,K(a_n,x)
\end{align*}
for all $x\in\R$. 
It remains to refer to Lemmas~\ref{lem:sqrt3>x>sqrt2} and \ref{lem:x>sqrt3}.
\end{proof}

Lemma~\ref{lem:sqrt3>x>sqrt2} is based on a series of other lemmas. To state those lemmas,
more notation is needed:
\begin{equation}
x_*:=\sqrt{\frac{5+2\sqrt{6}}{9-2\sqrt{6}}}
=\sqrt{\frac{23+28\sqrt{\frac23}}{19}}
\approx1.55;\label{eq:xx} 	
\end{equation}
\begin{align*} 
	\overline{R}&:=\{(a,x)\in\R^2\colon 0\le a\le1,\sqrt2\le x\le\sqrt3\};\\		
	R&:=\{(a,x)\in\R^2\colon 0<a<1,\sqrt2<x<\sqrt3\};\\	\LLe&:=\{(a,x)\in R\colon u<\sqrt2,v\le\sqrt3\};\\
\LG&:=\{(a,x)\in R\colon u<\sqrt2,v>\sqrt3\};\\
\GL_1&:=\{(a,x)\in R\colon u>\sqrt2,v<\sqrt3,x\le x_*\};\\
\GL_2&:=\{(a,x)\in R\colon u>\sqrt2,v<\sqrt3,x>x_*\};\\
\GG1&:=\{(a,x)\in R\colon u>\sqrt2,v>\sqrt3,a<\tfrac1{\sqrt3}\};\\
\GG2&:=\{(a,x)\in R\colon u>\sqrt2,v>\sqrt3,a\ge\tfrac1{\sqrt3}\};\\
\GE&:=\{(a,x)\in R\colon u>\sqrt2,v=\sqrt3\};\\
\ELe&:=\{(a,x)\in R\colon u=\sqrt2,v\le\sqrt3\};\\
\EG1&:=\{(a,x)\in R\colon u=\sqrt2,v>\sqrt3,a<\tfrac1{\sqrt3}\};\\
\EG2&:=\{(a,x)\in R\colon u=\sqrt2,v>\sqrt3,a\ge\tfrac1{\sqrt3}\};\\
\A1&:=\{(a,x)\in\overline{R}\colon a=0\};\\
\A2&:=\{(a,x)\in\overline{R}\colon a=1\};\\
\X{1,1}&:=\{(a,x)\in\overline{R}\colon x=\sqrt2,0<a<\tfrac1{\sqrt2}\};\\
\X{1,2}&:=\{(a,x)\in\overline{R}\colon x=\sqrt2,\tfrac1{\sqrt2}\le a<\tfrac{2\sqrt2}3\};\\
\X{1,3}&:=\{(a,x)\in\overline{R}\colon x=\sqrt2, \tfrac{2\sqrt2}3\le a<1\};\\
\X2&:=\{(a,x)\in\overline{R}\colon x=\sqrt3\},
\end{align*}
where $u$ and $v$ are defined by \eqref{eq:u} and \eqref{eq:v}.
Here, for example, the L in the first position in symbol $\LLe$ refers to ``less than'' in inequality $u<\sqrt2$, while the ligature $\mathrm{L\!e}$ in the second position refers to ``less than or equal to'' in inequality $v\le\sqrt3$.
Similarly, G and E in this notation refer to ``greater than'' and ``equal to'', respectively.
Symbol A refers here to a fixed value of $a$, and X to a fixed value of $x$.

It will be understood that the function $K$ is extended to $\A2$ by continuity, so that 
\begin{equation}\label{eq:K on A2}
K(a,x):=-2g(x)\quad\forall(a,x)\in\A2.	
\end{equation}

\begin{lemma}[$\LLe$]\label{lem:lle}
The function $K$ does not attain a maximum on $\LLe$.
\end{lemma}

\begin{lemma}[$\LG$]\label{lem:lg}
The function $K$ does not attain a maximum on $\LG$.
\end{lemma}

\begin{lemma}[$\GL_1$]\label{lem:gl1}
The function $K$ does not attain a maximum on $\GL_1$.
\end{lemma}

\begin{lemma}[$\GL_2$]\label{lem:gl2}
The function $K$ does not attain a maximum on $\GL_2$.
\end{lemma}

\begin{lemma}[$\GG1$]\label{lem:gg1}
The function $K$ does not attain a maximum on $\GG1$.
\end{lemma}

\begin{lemma}[$\GG2$]\label{lem:gg2}
The function $K$ does not attain a maximum on $\GG2$.
\end{lemma}

\begin{lemma}[$\GE$]\label{lem:ge}
One has $K\le0$ on $\GE$.
\end{lemma}

\begin{lemma}[$\ELe$]\label{lem:ele}
One has $K\le0$ on $\ELe$.
\end{lemma}

\begin{lemma}[$\EG1$]\label{lem:eg1}
One has $K\le0$ on $\EG1$.
\end{lemma}

\begin{lemma}[$\EG2$]\label{lem:eg2}
One has $K\le0$ on $\EG2$.
\end{lemma}

\begin{lemma}[$\A1$]\label{lem:a1}
One has $K=0$ on $\A1$.
\end{lemma}

\begin{lemma}[$\A2$]\label{lem:a2}
One has $K\le0$ on $\A2$.
\end{lemma}

\begin{lemma}[$\X{1,1}$]\label{lem:x11}
One has $K\le0$ on $\X{1,1}$.
\end{lemma}

\begin{lemma}[$\X{1,2}$]\label{lem:x12}
One has $K\le0$ on $\X{1,2}$.
\end{lemma}

\begin{lemma}[$\X{1,3}$]\label{lem:x13}
One has $K\le0$ on $\X{1,3}$.
\end{lemma}

\begin{lemma}[$\X2$]\label{lem:x2}
One has $K\le0$ on $\X2$.
\end{lemma}

The proofs of some of these lemmas require a great amount of symbolic and numerical computation. We have done that with the help of
Mathematica\texttrademark\ 5.2, which is rather effective and allows complete and easy control over the accuracy.
The idea behind these calculations is to reduce a problem involving transcendental inequalities and/or equations to 
an algebraic setting, for 
it is well known that the existence of a real-valued solution to any system of algebraic equations and/or inequalities can be determined in a completely algorithmic manner, according to a result by Tarski \cite{tarski,loja,coll}. 
This algorithm is implemented in Mathematica 5.2 via \texttt{Reduce} and related commands. 
However, execution of such commands may take a long time if the algebraic system is more than a little complicated; in such cases, Mathematica can use some human help.

\subsection{Proofs of the lemmas}\label{proofs of lemmas}

\begin{proof}[Proof of Lemma~\ref{lem:c1}]
Let $n=2$ and $a_1=a_2=\frac1{\sqrt2}$. Then
$\P(S_n\ge\sqrt2\,)=\frac14=c_1\P(Z\ge\sqrt2\,).$
\end{proof}

\begin{proof}[Proof of Lemma~\ref{lem:g?h}]
This follows from the well-known and easy-to-prove fact that the inverse Mills ratio $r$ is increasing.
\end{proof}

\begin{proof}[Proof of Lemma~\ref{lem:h1<h}]
On interval $(-\infty,0]$, one has
$h_1=1<1.6<h(0)\le h$, since $h$ is decreasing. 

On interval $(0,1]$, one similarly has
$h_1=\frac12<0.51<h(1)\le h$. 

On interval $[\sqrt2,\sqrt3\,]$, one has
$h_1=g\le h$, by Lemma~\ref{lem:g?h}. 

On interval $[\sqrt3,\infty)$, one has
$h_1=h$.

It remains to consider the interval $(1,\sqrt2)$. 
For
$x\in(1,\sqrt2)$, one has $h_1(x)=p(x):=\frac1{2x^2}$. 
Now let us apply a l'Hospital-type rule for monotonicity; see e.g.\ \cite[Proposition 4.3]{borwein} and the notation used therein.
One has $p(\infty-)=h(\infty-)=0$ and, for some constant $C>0$,
$h'(x)/p'(x)=Cx^3\vpi(x),$ so that $\frac{h'}{p'}\u\d$ on $(0,\infty)$. Hence, $\frac hp\d$ or $\u\d$ on $(0,\infty)$, and so, 
$\inf_{(1,\sqrt2)}\frac h{h_1}
=\min_{[1,\sqrt2]}\frac hp
=\min_{\{1,\sqrt2\}}\frac hp\approx1.01>1$,
whence $h_1<h$ on $(1,\sqrt2)$.  
\end{proof}

\begin{proof}[Proof of Lemma~\ref{lem:sqrt3>x>sqrt2}]
Since the function $h_1$ is upper-semicontinuous on $\R$ and continuous on $[\sqrt2,\infty)$, the function $K$ is upper-semicontinuous on the compact rectangle $\overline{R}$ and therefore attains its maximum on $\overline{R}$.
Now Lemma~\ref{lem:sqrt3>x>sqrt2} follows from Lemmas~\ref{lem:lle}--\ref{lem:x2} (proved below), since $\overline{R}$ is the union of sets $\LLe,\dots,\X2$.
\end{proof}

\begin{proof}[Proof of Lemma~\ref{lem:x>sqrt3}]
Let $0\le a<1$ and $x\ge\sqrt3$.
Then $u\ge\sqrt2$ and $v\ge\sqrt3$, where $u$ and $v$ are defined by \eqref{eq:u} and \eqref{eq:v}, as before. Therefore and in view of Lemma~\ref{lem:h1<h} and \eqref{eq:h1}, one has $h_1(u)\le h(u)$, $h_1(v)=h(v)$,  and $h_1(x)=h(x)$, so that
$$K(a,x)\le2c_2\cdot\big(\tfrac12\ps(u)+\tfrac12\ps(v)-\ps(x)\big)\le0,$$
as shown in the mentioned proof in \cite{bgh}. 
\end{proof}

\begin{proof}[Proof of Lemma~\ref{lem:lle} ($\LLe$)]
Expressing $x$ and $u$ in view of \eqref{eq:u} and \eqref{eq:v} in terms of $a$ and $v$ as 
$x(a,v):=\sqrt{1-a^2}\,v-a$ and 
$\tilde u(a,v):=v-\frac{2a}{\sqrt{1-a^2}}$, respectively, one has
\begin{equation}
	K(a,x)=k(a,v):=k_{\LLe}(a,v):= 
\frac1{2\tilde u(a,v)^2}+g(v)-2g(x(a,v))
\quad\forall(a,x)\in\LLe,
\label{eq:kLLe}
\end{equation}
since $u>1$ and $v>\sqrt2$ on $R$.
For $(a,x)\in R$, let
\begin{align}
(D_a k)(a,v)&:=4\ps(\sqrt2)\big(x(a,v)-a\big)^3\,\p k a (a,v); \label{eq:dak}\\
(D_{a,a}k)(a,x)&:=
\frac{125(1-a^2)^2}{(x-a)^2\vpi(x)}
\p{(D_a k)}a\big(a,v(a,x)\big).\label{eq:daak}
\end{align}
Then $(D_{a,a}k)(a,x)$ is an algebraic expression (and even a polynomial) in $a$ and $x$. 
With Mathematica 5.2, one can therefore use the command
$${\bf\texttt{Reduce[daakx<=0 \&\& 0<a<1 \&\& Sqrt[2]<x<Sqrt[3]]}}$$ 
(where {\bf\texttt{daakx}} 
stands for $(D_{a,a}k)(a,x)$), which outputs \verb2False2, meaning that $D_{a,a}k>0$ on $R$.  
By \eqref{eq:daak}, this implies $\p{(D_a k)}a\big(a,v(a,x)\big)>0$ for $(a,x)\in R$, so that $(D_a k)(a,v)$ is increasing in $a$ for every fixed value of $v$; more exactly, $(D_a k)(a,v)$ is increasing in $a\in\big(a_1(v),a_2(v)\big)$ for every fixed value of $v\in(\sqrt2,\sqrt3]$, where $a_1$ and $a_2$ are certain functions, such that 
$$\big(a,x(a,v)\big)\in\LLe \iff
\Big(v\in(\sqrt2,\sqrt3]\ \&\ a\in\big(a_1(v),a_2(v)\big)\Big).$$
Thus, for every fixed value of $v\in(\sqrt2,\sqrt3]$ the sign pattern of $(D_a k)(a,v)$ in $a\in\big(a_1(v),a_2(v)\big)$ is $-$ or $+$ or $-+$; that is, $(D_a k)(a,v)$ may change sign only from $-$ to $+$ as $a$ increases. By \eqref{eq:dak}, $\p k a (a,v)$ has the same sign pattern.
Hence, for every fixed value of $v\in(\sqrt2,\sqrt3]$ one has
$k(a,v)\d$ or $\u$ or $\d\u$ in $a\in\big(a_1(v),a_2(v)\big)$. 
Now Lemma~\ref{lem:lle} follows. 
\end{proof}

\begin{proof}[Proof of Lemma~\ref{lem:lg} ($\LG$)]
This proof is almost identical to that of Lemma~\ref{lem:lle}, except that the term $g(v)$ in \eqref{eq:kLLe} is replaced here by $h(v)$, and the interval $(\sqrt2,\sqrt3]$ is replaced by $(\sqrt3,5\sqrt2\,)$. However, both terms $g(v)$ and $h(v)$ are constant for any fixed value of $v$.
\end{proof}

\begin{proof}[Proof of Lemma~\ref{lem:gl1} ($\GL_1$)]
One has
\begin{equation}\label{eq:GL1}
(a,x)\in\GL_1 \iff
\Big(\sqrt2<x\le x_*\ \&\ 0<a<a_1(x)\Big),
\end{equation}
where $a_1(x):=\frac x3-\frac13\sqrt{6-2x^2}$.
Since $\sqrt2<u<v<\sqrt3$ on $\GL_1$ (where again $u$ and $v$ are defined by \eqref{eq:u} and \eqref{eq:v}), one has 
\begin{equation}\label{eq:kGL1}
	K(a,x)=k(a,x):=k_{\GL}(a,x):= 
g(u)+g(v)-2g(x)
\quad\forall(a,x)\in\GL_1.
\end{equation}
For $(a,x)\in R$, let
\begin{align}
(D_a k)(a,x)&:=\p k a (a,x)\cdot
\frac{125\,\ps(\sqrt2\,)\left(1-a^2\right)^{3/2}}
{(1+a x)(251+v) \vpi(v)};\label{eq:dakGL1}\\
(D_{a,a}k)(s,x)&:=
\p{(D_a k)}a(a,x)\cdot
(1-a^2)^3 (1+ax)^2 (251+v)^2 e^{-\frac{2ax}{1-a^2}}, \label{eq:daakGL1}
\end{align}
where $\sqrt{1-s^2}$ is substituted for $a$ in the right-hand side of \eqref{eq:daakGL1}.
Then $(D_{a,a}k)(s,x)$ is a polynomial in $s$ and $x$. 
Using again the Mathematica command 
\texttt{Reduce}, namely
\begin{equation}\label{eq:ReddaakGL1}
	\texttt{Reduce[daak[s,x]>0 \&\& Sqrt[2]<x<Sqrt[3] \&\& 0<s<1]},
\end{equation}
where \texttt{daak[s,x]} stands for $(D_{a,a}k)(s,x)$),
one sees that for every $(s,x)\in R$
\begin{equation}\label{eq:ss1,ss2}
(D_{a,a}k)(s,x)>0 \iff 
\Big(\sqrt2<x<x_{**}\ \&\ s_{*,1}(x)<s<s_{*,2}(x)\Big),	
\end{equation}
where $x_{**}$ is a certain number between $\sqrt2$ and $\sqrt3$, and $s_{*,1}$ and $s_{*,2}$ are certain functions. 
\big(In fact, $x_{**}\approx1.678696$ 
is a root of a certain polynomial of degree $32$ and, for each $x\in(\sqrt2,x_{**})$, the values $s_{*,1}(x)$ and $s_{*,2}(x)$ are two of the roots $s$ of the polynomial $(D_{a,a}k)(s,x)$.\big) 
Next, \texttt{Reduce[daak[95/100,x]<=0 \&\& Sqrt[2]<x<=xx]} produces \texttt{False}; here \texttt{xx} stands for $x_*$ -- recall the definition of $\GL_1$ and \eqref{eq:xx}; that is,
$(D_{a,a}k)(\frac{95}{100},x)>0$ $\forall x\in(\sqrt2,x_*]$.
Hence and in view of \eqref{eq:ss1,ss2}, 
\begin{equation}\label{eq:95/100}
s_{*,1}(x)<\tfrac{95}{100}<s_{*,2}(x)\quad\forall x\in(\sqrt2,x_*].	
\end{equation}
On the other hand, setting 
\begin{equation}\label{eq:s1}
	s_1(x):=\sqrt{1-a_1(x)^2}
\end{equation}
with $a_1$ as in \eqref{eq:GL1}, one has $s_1>\tfrac{95}{100}$ on $(\sqrt2,x_*]$. 
Hence, by \eqref{eq:95/100}, $s_1>s_{*,1}$ on $(\sqrt2,x_*]$. 
So, in view of \eqref{eq:ss1,ss2}, the sign pattern of $(D_{a,a}k)(s,x)$ in $s\in(s_1(x),1)$ is $-$ or $+-$, depending on whether $s_1(x)\ge s_{*,2}(x)$ or not, for each $x\in(\sqrt2,x_*]$. 
Now \eqref{eq:daakGL1} and \eqref{eq:s1} imply that 
the sign pattern of $\p{(D_a k)}a(a,x)$ in $a\in(0,a_1(x))$ is $-$ or $-+$, so that
$(D_a k)(a,x)\d$ or $\d\u$ in $a\in(0,a_1(x))$, for each $x\in(\sqrt2,x_*]$.
Also, $(D_a k)(0,x)=0$ for all $x\in\R$. 
So, the sign pattern of $(D_a k)(a,x)$ in $a\in(0,a_1(x))$ is $-$ or $-+$, for each $x\in(\sqrt2,x_*]$; 
in view of \eqref{eq:dakGL1}, $\p k a(a,x)$ has the same sign pattern.
Thus, $k(a,x)\d$ or $\d\u$ in $a\in(0,a_1(x))$, for each $x\in(\sqrt2,x_*]$.
Recalling \eqref{eq:GL1}, we complete the proof of Lemma~\ref{lem:gl1}.
\end{proof}

\begin{proof}[Proof of Lemma~\ref{lem:gl2} ($\GL_2$)]
This proof is similar to that of Lemma~\ref{lem:gl1}.
Here one has
\begin{equation}\label{eq:GL2}
(a,x)\in\GL_2 \iff
\Big(x_*<x<\sqrt3\ \&\ 0<a<a_2(x)\Big),
\end{equation}
where $a_2(x):=-\frac x4+\frac14\sqrt{12-3x^2}$.
Relation \eqref{eq:kGL1} holds here, and we retain definitions \eqref{eq:dakGL1} and \eqref{eq:daakGL1}.
Definition \eqref{eq:s1} is replaced here by
\begin{equation}\label{eq:s2}
	s_2(x):=\sqrt{1-a_2(x)^2}.
\end{equation}
Letting 
$$s_*:=s_2(x_*)=\frac1{12}\sqrt{\frac{1728+384\sqrt{6}}{19}}\approx0.98761,$$
one can see that
\begin{equation}\label{eq:s2>ss}
	s_2(x)\ge s_*\quad\forall x\in(x_*,\sqrt3\,).
\end{equation}
On the other hand, using instead of \eqref{eq:ReddaakGL1}
the Mathematica command
\begin{equation*}
	\texttt{Reduce[daak[s,x]>0 \&\& xx<x<Sqrt[3] \&\& ss<s<1, Quartics->True]}, 
\end{equation*}
where $\texttt{xx}$ stands for $x_*$ and $\texttt{ss}$ stands for $s_*$, 
one sees that 
\begin{multline}\label{eq:ss,ss2}
\Big((D_{a,a}k)(s,x)>0
\ \&\ x_*<x<\sqrt3\ \&\ s_*<s<1\Big) \\
\iff 
\Big(x_*<x<x_{***}\ \&\ s_*<s<s_{*,2}(x)\Big),	
\end{multline}
where $x_{***}\approx1.678694$ is a root of a certain polynomial of degree $20$ and
$s_{*,2}(x)$ is the same root in $s$ of the polynomial $(D_{a,a}k)(s,x)$ as $s_{*,2}(x)$ in \eqref{eq:ss1,ss2}.   
Hence, in view of \eqref{eq:ss,ss2} and \eqref{eq:s2>ss}, the sign pattern of $(D_{a,a}k)(s,x)$ in $s\in(s_2(x),1)$ is $-$ or $+-$, for each $x\in(x_*,\sqrt3\,)$. 

Now \eqref{eq:daakGL1} and \eqref{eq:s2} imply that 
the sign pattern of $\p{(D_a k)}a(a,x)$ in $a\in(0,a_2(x))$ is $-$ or $-+$, so that
$(D_a k)(a,x)\d$ or $\d\u$ in $a\in(0,a_2(x))$, for each $x\in(x_*,\sqrt3\,)$.
Also, $(D_a k)(0,x)=0$ for all $x\in\R$. 
So, the sign pattern of $(D_a k)(a,x)$ in $a\in(0,a_2(x))$ is $-$ or $-+$, for each $x\in(x_*,\sqrt3\,)$; 
in view of \eqref{eq:dakGL1}, $\p k a(a,x)$ has the same sign pattern.
Thus, $k(a,x)\d$ or $\d\u$ in $a\in(0,a_2(x))$, for each $x\in(x_*,\sqrt3\,)$.
Recalling \eqref{eq:GL2}, we complete the proof of Lemma~\ref{lem:gl2}.
\end{proof}

\begin{proof}[Proof of Lemma~\ref{lem:gg1} ($\GG1$)]
Expressing $x$ and $v$ in view of \eqref{eq:u} and \eqref{eq:v} in terms of $a$ and $u$ as 
$\tilde x(a,u):=\sqrt{1-a^2}\,u+a$ and 
$\tilde v(a,u):=u+\frac{2a}{\sqrt{1-a^2}}$, respectively, 
and taking into account that $u<\sqrt3$ on $\GG1$,
one has
\begin{equation}
	K(a,x)=k(a,u):=k_{\GG1}(a,u):= 
g(u)+h(\tilde v(a,u))-2g(\tilde x(a,u))
\quad\forall(a,x)\in\GG1.
\label{eq:kGG1}
\end{equation}
For $(a,u)\in R$, let
\begin{align}
d(a,u)&:=
\p ka(a,u)\cdot125\,\ps(\sqrt2)\,(1-a^2)^{3/2}
/\vpi(\tilde v(a,u));
   \label{eq:dakGG1}\\
d_a(a,u)&:=
\p da(a,u)\cdot(1-a^2)\,\frac{\vpi(\tilde v(a,u))}{\vpi(\tilde x(a,u))};\notag
\\
d_u(a,u)&:=
\p du(a,u)\cdot
\frac{\vpi(\tilde v(a,u))}{\vpi(\tilde x(a,u))\sqrt{1-a^2}}. \notag
\end{align}
Then $d_a(a,u)$ and $d_u(a,u)$ are polynomials in $a$, $u$, and $\sqrt{1-a^2}$. Using 
\begin{multline*}
\texttt{Reduce[da[a,u]==0 \&\& du[a,u]==0 \&\& 0<a<1/Sqrt[3] \&\&} \\ 
\texttt{Sqrt[2]<u<Sqrt[3], \{a,u\}, Reals]},
\end{multline*}
where $\texttt{da[a,u]}$ and $\texttt{du[a,u]}$ stand for $d_a(a,u)$ and $d_u(a,u)$, one sees that the system of equations $d_a(a,u)=0=d_u(a,u)$ has a unique solution $(a_*,u_*)\approx(0.11918,1.57770)$ in $(a,u)\in(0,\frac1{\sqrt3})\times(\sqrt2,\sqrt3\,)$, and $d(a_*,u_*)\approx-0.44<0$.

Let us consider next the values of $d$ on the boundary of the rectangle $(0,\frac1{\sqrt3})\times(\sqrt2,\sqrt3\,)$.

First, $d(0,u)$ is an increasing affine function of $u$, and $d(0,\sqrt3\,)\approx-0.4<0$. Hence, $d(0,u)<0$ $\forall u\in(\sqrt2,\sqrt3\,)$.

Second,
\begin{multline*}
\frac{27}{2\,e^{(5+4\sqrt2u+u^2)/6}}\,
\p du(\tfrac1{\sqrt3},u) \\
=-\sqrt{2}u^3-\left(3+251 \sqrt{3}\right) u^2-\sqrt{2}
   \left(3+251 \sqrt{3}\right) u+251 \sqrt{3}+7
\end{multline*}
is decreasing in $u$ and takes on value $-9-753\sqrt{3}<0$ at $u=\sqrt2$, so that $\p du(\tfrac1{\sqrt3},u)<0$ for $u\in(\sqrt2,\sqrt3\,)$
and $d(\tfrac1{\sqrt3},u)\d$ in $u\in(\sqrt2,\sqrt3\,)$. Moreover, $d(\tfrac1{\sqrt3},\sqrt2)<0$. Thus, $d(\tfrac1{\sqrt3},u)<0$ $\forall u\in(\sqrt2,\sqrt3\,)$. 

Third,
$$
\left(1-a^2\right)\,
\frac{\vpi(\tilde v(a,\sqrt2\,))}{\vpi(\tilde x(a,\sqrt2\,))}\,
\p da(a,\sqrt2) 
=p_1(a)+\sqrt{1-a^2}\,p_2(a),
$$
where $p_1(a)$ and $p_2(a)$ are certain polynomials in $a$.
Therefore, the roots $a$ of $\p da(a,\sqrt2)$ are among the roots of
the polynomial $p_{1,2}(a):=p_1(a)^2-(1-a^2)p_2(a)^2$, which has exactly two roots $a\in(0,\tfrac1{\sqrt3})$. Of the latter roots, one is not a root of $\p da(a,\sqrt2)$. 
Also, $\p da(0,\sqrt2)=1>0$ and $\p da(\tfrac1{\sqrt3},\sqrt2)
<0$.
Hence, $\p da(a,\sqrt2)$ has exactly one root, $a_*\approx0.2224$, in $a\in(0,\tfrac1{\sqrt3})$ and, moreover, $d(a,\sqrt2)\u$ in $a\in(0,a_*]$ and $d(a,\sqrt2)\d$ in $a\in[a_*,\tfrac1{\sqrt3})$. 
Besides, $d(a_*,\sqrt2)\approx-0.088<0$. 
Thus, $d(a,\sqrt2)<0$ $\forall a\in(0,\tfrac1{\sqrt3})$. 

Fourth (very similar to third),
$$
\left(1-a^2\right)\,
\frac{\vpi(\tilde v(a,\sqrt3\,))}{\vpi(\tilde x(a,\sqrt3\,))}\,
\p da(a,\sqrt3\,) 
=p_1(a)+\sqrt{1-a^2}\,p_2(a),
$$
where $p_1(a)$ and $p_2(a)$ are certain polynomials in $a$, different from the polynomials $p_1(a)$ and $p_2(a)$ in the previous paragraph.
Therefore, the roots $a$ of $\p da(a,\sqrt3)$ are among the roots of
the polynomial $p_{1,2}(a):=p_1(a)^2-(1-a^2)p_2(a)^2$, which has exactly two roots $a\in(0,\tfrac1{\sqrt3})$. Of the latter roots, one is not a root of $\p da(a,\sqrt3\,)$. 
Also, $\p da(0,\sqrt3\,)=1>0$ and $\p da(\tfrac1{\sqrt3},\sqrt3\,)
<0$.
Hence, $\p da(a,\sqrt3\,)$ has exactly one root, $a_*\approx0.06651$, in $a\in(0,\tfrac1{\sqrt3})$ and, moreover, $d(a,\sqrt3\,)\u$ in $a\in(0,a_*]$ and $d(a,\sqrt3\,)\d$ in $a\in[a_*,\tfrac1{\sqrt3})$. 
Besides, $d(a_*,\sqrt3\,)\approx-0.358<0$. 
Thus, $d(a,\sqrt3\,)<0$ $\forall a\in(0,\tfrac1{\sqrt3})$. 

We conclude that $d(a,u)<0$  $\forall(a,u)\in[0,\frac1{\sqrt3}]\times[\sqrt2,\sqrt3\,]$.
By \eqref{eq:dakGG1}, the same holds for $\p ka(a,u)$.
It remains to recall \eqref{eq:kGG1} and note that 
$$(a,x)\in\GG1\implies
\big(a,u(a,x)\big)
\in(0,\tfrac1{\sqrt3})\times(\sqrt2,\sqrt3\,).$$ 
\end{proof}

\begin{proof}[Proof of Lemma~\ref{lem:gg2} ($\GG2$)]
In view of Lemma~\ref{lem:g?h},
\begin{equation}
K(a,x)= 
(g\wedge h)(u(a,x))+h(v(a,x))-2g(x)
\quad\forall(a,x)\in\GG2.
\label{eq:kGG2}
\end{equation}
One has $\p ua>0$ on $\GG2$ and $\p va>0$ on $R$.
Since $g=
\frac{c_1}{250}\,\big(251\ps+\vpi\big)$ and $h=c_2\,\ps$ are decreasing on $[0,\infty)$, we conclude that
$K(a,x)$ is decreasing in $a$ on $\GG2$.
\end{proof}

\begin{proof}[Proof of Lemma~\ref{lem:ge} ($\GE$)]
One has
\begin{equation*}
	(a,x)\in\GE \iff
	\Big( x_*<x<\sqrt3\ \&\  a=a_2(x):=\tfrac{1}{4} \sqrt{12-3 x^2}-\tfrac{x}{4}\;\Big),
\end{equation*}
where, as before, $x_*$ is defined by \eqref{eq:xx}.
Therefore, for all $(a,x)\in\GE$
\begin{equation*}
	K(a,x)=k(x):=k_{\GE}(x):=K(a_2(x),x)
	=g\big(u(a_2(x),x)\big)+g(\sqrt3)-2g(x), 
\end{equation*}
and it suffices to show that $k\le0$ on $[x_*,\sqrt3]$.
For $x\in[x_*,\sqrt3]$,
let
\begin{align*}
k_1(x)&:=\frac
{k'(x)\,\ps(\sqrt2\,)}
{\vpi(x)(251+x)};\\
k_2(x)&:=k'_1(x)\cdot
\frac{
125\,(x+251)^2\, \left(4-x^2\right)^{3/2} \rho(x)^{13/2}
   }
   {16\sqrt2\,\vpi\big(u(a_2(x),x)\big)/\vpi(x)},
\end{align*}
where $\rho(x):=x^2+x\sqrt{3} \sqrt{4-x^2}+2$.
Then $k_2(x)=p_1(x)+\sqrt{\rho(x)}\,p_2(x)$, where $p_1(x)$ and $p_2(x)$ are some polynomials in $x$ and $\sqrt{4-x^2}$.
Hence, the roots of $k_2(x)$ are among the roots of 
$$p_{1,2}(x):=p_1(x)^2-\rho(x)p_2(x)^2
=p_{1,2,1}(x)+\sqrt{4-x^2}\,p_{1,2,2}(x),$$
where $p_{1,2,1}(x)$ and $p_{1,2,2}(x)$ are some polynomials in $x$.
Hence, 
the roots of $k_2(x)$ are among the roots of 
$$
\frac{
p_{1,2,1}(x)^2-(4-x^2)p_{1,2,2}(x)^2
}
{1024\, \left(x^2-1\right)^{14}},$$
which is a polynomial of degree $32$ and has exactly one root 
in $[x_*,\sqrt3]$.
Also, $k_2(x_*)\approx1.39\times10^7>0$ and 
$k_2(\sqrt3\,)\approx-2.07\times10^6<0$.
Therefore, $k_2$ and hence $k'_1$ have the sign pattern $+-$ on $[x_*,\sqrt3]$.
Next, $k_1(x_*)\approx-4.8494<0$ and $k_1(\sqrt3\,)=0$, so that $k_1$ and hence $k'$ have the sign pattern $-+$ on $[x_*,\sqrt3]$.
It follows that $k$ does not have a local maximum on $(x_*,\sqrt3)$. 
At that, $k(x_*)\approx-3.0133\times10^{-6}<0$ and $k(\sqrt3\,)=0$.
Thus, $k\le0$ on $[x_*,\sqrt3]$.
\end{proof}

\begin{proof}[Proof of Lemma~\ref{lem:ele} ($\ELe$)]
This proof is very similar to that of Lemma~\ref{lem:ge} ($\GE$).
One has
\begin{equation}\label{eq:ELe}
	(a,x)\in\ELe \iff
	\Big( \sqrt2<x\le x_*\ \&\  a=a_1(x):=\tfrac{x}{3}-\tfrac{1}{3} \sqrt{6-2 x^2}\;\Big),
\end{equation}
where, as before, $x_*$ is defined by \eqref{eq:xx}.
Therefore, for all $(a,x)\in\LLe$
\begin{equation*}
	K(a,x)=k(x):=k_{\ELe}(x):=K(a_1(x),x)
	=g(\sqrt2)+g\big(v(a_1(x),x)\big)-2g(x), 
\end{equation*}
and it suffices to show that $k\le0$ on $[\sqrt2,x_*]$.
For $x\in[\sqrt2,x_*]$,
let
\begin{align*}
k_1(x)&:=\frac
{500\,\ps(\sqrt2)\,k'(x)}
{\vpi(x)(251+x)};\\
k_2(x)&:=k'_1(x)\cdot
\frac{
\sqrt{2} (x+251)^2 \left(3-x^2\right)^{3/2} \rho(x)^{13/2}
   }
   {9\vpi\big(v(a_1(x),x)\big)/\vpi(x)},
\end{align*}
where $\rho(x):=x^2+2x\sqrt{2} \sqrt{3-x^2}+3$.
Then $k_2(x)=p_1(x)+\sqrt{\rho(x)}\,p_2(x)$, where $p_1(x)$ and $p_2(x)$ are some polynomials in $x$ and $\sqrt{3-x^2}$.
Hence, the roots of $k_2(x)$ are among the roots of 
$$p_{1,2}(x):=p_1(x)^2-\rho(x)\,p_2(x)^2
=p_{1,2,1}(x)+\sqrt{3-x^2}\,p_{1,2,2}(x),$$
where $p_{1,2,1}(x)$ and $p_{1,2,2}(x)$ are some polynomials in $x$.
Hence, 
the roots of $k_2(x)$ are among the roots of 
$$
\frac{
p_{1,2,1}(x)^2-(3-x^2)p_{1,2,2}(x)^2
}
{125524238436 \left(x^2-1\right)^{14}},$$
which is a polynomial of degree $32$ and has exactly one root 
in $[\sqrt2,x_*]$.
Also, $k_2(\sqrt2)\approx-6.32\times10^7<0$ and 
$k_2(x_*)\approx1.06\times10^8>0$.
Therefore, $k_2$ and hence $k'_1$ have the sign pattern $-+$ on $[\sqrt2,x_*]$.
Next, $k_1(\sqrt2)=0$ and $k_1(x_*)\approx0.000426>0$, so that $k_1$ and hence $k'$ have the sign pattern $-+$ on $[\sqrt2,x_*]$.
It follows that $k$ does not have a local maximum on $(\sqrt2,x_*)$. 
At that, $k(\sqrt2)=0$ and $k(x_*)\approx-3.0133\times10^{-6}<0$.
Thus, $k\le0$ on $[\sqrt2,x_*]$.
\end{proof}

\begin{proof}[Proof of Lemma~\ref{lem:eg1} ($\EG1$)]
This proof is similar to that of Lemma~\ref{lem:ele}.
One has
\begin{equation}\label{eq:EG1}
	(a,x)\in\EG1 \iff
	\Big( x_*<x<\sqrt3\ \&\  a=a_1(x):=\tfrac{x}{3}-\tfrac{1}{3} \sqrt{6-2 x^2}\;\Big).
\end{equation}
Therefore, for all $(a,x)\in\EG1$
\begin{equation*}
	K(a,x)=k(x):=k_{\EG1}(x):=K(a_1(x),x)
	=g(\sqrt2\,)+h\big(v(a_1(x),x)\big)-2g(x), 
\end{equation*}
and it suffices to show that $k\le0$ on $[x_*,\sqrt3]$.
For $x\in[x_*,\sqrt3]$,
let
\begin{align*}
k_1(x)&:=\frac
{500\,\ps(\sqrt2)\,k'(x)}
{\vpi(x)(251+x)};\\
k_2(x)&:=k'_1(x)\cdot
\frac{\sqrt2\,(x+251)^2 \left(3-x^2\right)^{3/2} 
\rho(x)^{9/2}  
 }{9\,r(\sqrt{3}\,)
\vpi\big(v(a_1(x),x)\big)/\vpi(x)
   },
\end{align*}
where $\rho(x):=x^2+2x\sqrt{2} \sqrt{3-x^2}+3$.

Then $k_2(x)=p_1(x)+\sqrt{3-x^2}\,p_2(x)$, where $p_1(x)$ and $p_2(x)$ are some polynomials in $x$.
Hence, the roots of $k_2(x)$ are among the roots of 
$$p_{1,2}(x):=
\frac{p_1(x)^2-(3-x^2)p_2(x)^2}
{4374\,\big(1+251/r(\sqrt3\,)\big)^2\,(x^2-1)^4}
,$$
which is a polynomial of degree $14$ and has exactly one root 
in $[x_*,\sqrt3]$, $x_{\#}\approx1.6012$.
Also,  
$k_2(x_*)\approx1.1722\times10^6>0$
and $k_2(\sqrt3)\approx-3.8778\times10^7<0$.
Therefore, $k_2$ and hence $k'_1$ have the sign pattern $+-$ on $[x_*,\sqrt3]$, so that
$\max_{[x_*,\sqrt3]}k_1=k_1(x_{\#})\approx-0.00034907<0$.
It follows that 
$k'<0$ and hence $k\d$ on $[x_*,\sqrt3]$.
At that, $k(x_*)\approx-3.0133\times10^{-6}<0$.
Thus, $k\le0$ on $[x_*,\sqrt3]$.
\end{proof}

\begin{proof}[Proof of Lemma~\ref{lem:eg2} ($\EG2$)]
One has
\begin{equation}\label{eq:EG2}
	(a,x)\in\EG2 \iff
	\Big( \sqrt2<x<\sqrt3\ \&\  a=a_2(x):=\tfrac{x}{3}+\tfrac{1}{3} \sqrt{6-2 x^2}\;\Big).
\end{equation}
Therefore, for all $(a,x)\in\EG2$
\begin{equation*}
	K(a,x)=k(x):=k_{\EG2}(x):=K(a_2(x),x)
	=g(\sqrt2)+h\big(v(a_2(x),x)\big)-2g(x), 
\end{equation*}
and it suffices to show that $k\le0$ on $(\sqrt2,\sqrt3)$.

For the functions $a_1$ (defined in \eqref{eq:ELe} and \eqref{eq:EG1}) and $a_2$ (defined in \eqref{eq:EG2}), and for
$x\in(\sqrt2,\sqrt3)$, one has $a_2(x)\ge a_1(x)$; also, $h(z)$ is decreasing in $z$ and $v(a,x)$ is increasing in $a$. Hence, $h\big(v(a_2(x),x)\big)\le h\big(v(a_1(x),x)\big)$, so that $k_{\EG2}\le k_{\EG1}$ on $[x_*,\sqrt3)$.

Similarly, in view of Lemma~\ref{lem:g?h} one has 
$h\big(v(a_2(x),x)\big)\le g\big(v(a_2(x),x)\big)\le g\big(v(a_1(x),x)\big)$ $\forall x\in(\sqrt2,x_*]$, so that
$k_{\EG2}\le k_{\ELe}$ on $(\sqrt2,x_*]$.

Now Lemma~\ref{lem:eg2} follows, because it was shown in the proofs of Lemmas~\ref{lem:ele} and \ref{lem:eg1}, respectively, that $k_{\ELe}\le0$ on $[\sqrt2,x_*]$ and $k_{\EG1}\le0$ on $[x_*,\sqrt3\,]$.
\end{proof}

\begin{proof}[Proof of Lemma~\ref{lem:a1} ($\A1$)]
This is trivial.
\end{proof}

\begin{proof}[Proof of Lemma~\ref{lem:a2} ($\A2$)]
This is also trivial, in view of \eqref{eq:K on A2}.
\end{proof}

\begin{proof}[Proof of Lemma~\ref{lem:x11} ($\X{1,1}$)]
On $\X{1,1}$, one has $u<\sqrt2\le v$. 
Also, by Lemma~\ref{lem:g?h}, $g\ge h$ on $[\sqrt3,\infty)$.
Therefore, for all $(a,x)\in\X{1,1}$
\begin{equation}\label{eq:kx11}
	K(a,x)\le k(a):=k_{\X{1,1}}(a):=
	\frac1{2u(a,\sqrt2\,)^2}+g\big(v(a,\sqrt2\,)\big)-2g(\sqrt2\,), 
\end{equation}
and it suffices to show that $k\le0$ on $[0,\frac1{\sqrt2})$.
For $a\in[0,\frac1{\sqrt2})$,
let
\begin{align*}
k_1(a)&:=\frac
{2000\,\ps(\sqrt2)\,k'(a)}
{\la(a)},\quad
\la(a):=\frac{1-a\sqrt{2}}{\left(\sqrt{2}-a\right)^3}>0;\\
k_2(a)&:=k'_1(a)\cdot
(\sqrt{2}-a)^4 \left(1-a^2\right)^4\,\la(a)^2
/\Big(\sqrt2\,\vpi\big(v(a,\sqrt2\,)\big)\Big).
\end{align*}
Then $k_2(a)=\sqrt{1-a^2}\,p_1(a)+p_2(a)$, where $p_1(a)$ and $p_2(a)$ are some polynomials in $a$.
Hence, the roots of $k_2(a)$ are among the roots of 
$$p_{1,2}(a):=
(1-a^2)p_1(a)^2-p_2(a)^2
,$$
which is a polynomial of degree $12$ and has exactly two roots 
in $[0,\frac1{\sqrt2})$.
Of those two roots, one is not a root of $k_2(a)$, so that $k_2(a)$
has at most one root in $[0,\frac1{\sqrt2})$.
Also,  
$k_2(0)=251\,\sqrt2>0$ and 
$k_2(\frac1{\sqrt2})=-127<0$.
Therefore, $k_2$ and hence $k'_1$ have the sign pattern $+-$ on $[0,\frac1{\sqrt2}]$, so that
$k_1\u\d$ on $[0,\frac1{\sqrt2}]$.
At that the values of $k_1$ at points $0$, $\frac6{10}$, and $\frac7{10}$ are approximately $-52<0$, $48>0$, and $-344<0$, respectively.
Therefore, $k_1$ and hence $k'$ have the sign pattern $-+-$ on $[0,\frac1{\sqrt2})$ and, moreover, 
the only local maximum of $k$ on $(0,\frac1{\sqrt2}]$ occurs only between $\frac6{10}$ and $\frac7{10}$;
in fact, it occurs at $a\approx0.67433$ and equals $\approx-0.00013578<0$. 
It remains to note that $k(0)\approx-0.0028660<0$.
\end{proof}

\begin{proof}[Proof of Lemma~\ref{lem:x12} ($\X{1,2}$)]
This proof is similar to that of Lemma~\ref{lem:gg2}. In place of \eqref{eq:kGG2}, here one still has relation \eqref{eq:kx11}
for all $(a,x)\in\X{1,2}$, since $u<\sqrt2\le v$ on
$\X{1,2}$ as well. 
Since $u(a,\sqrt2\,)\u$ and $v(a,\sqrt2\,)\u$
in $a\in[\frac1{\sqrt2},\frac{2\sqrt2}3)$, 
one has $k\d$ on $[\frac1{\sqrt2},\frac{2\sqrt2}3)$, so that
the maximum of $k$ on $[\frac1{\sqrt2},\frac{2\sqrt2}3)$ equals $k(\frac1{\sqrt2})$, which is negative, in view of Lemma~\ref{lem:x11} and the continuity of $k$.
\end{proof}

\begin{proof}[Proof of Lemma~\ref{lem:x13} ($\X{1,3}$)]
This proof is similar to that of Lemma~\ref{lem:x12}. In place of \eqref{eq:kx11}, here one has 
\begin{equation*}\label{eq:kx13}
	K(a,x)\le k(a):=k_{\X{1,3}}(a):=
	g\big(u(a,\sqrt2\,)\big)+g\big(v(a,\sqrt2\,)\big)-2g(\sqrt2\,), 
\end{equation*}
for all $(a,x)\in\X{1,3}$, since $u\ge\sqrt2$ and $v\ge\sqrt2$ on
$\X{1,3}$. 
Since $k\d$ in $a\in[\frac{2\sqrt2}3,1)$, the maximum of $k$ on $[\frac{2\sqrt2}3,1)$ equals $k(\frac{2\sqrt2}3)\approx-0.25287<0$.
\end{proof}

\begin{proof}[Proof of Lemma~\ref{lem:x2} ($\X2$)]
This follows immediately from Lemma~\ref{lem:x>sqrt3}. 
\end{proof}


\end{document}